\input amstex
\documentstyle{amsppt}
\magnification=1200
\hsize=5.2in
\vsize=7.5in

\topmatter
\title
 Estimates on Lattice Points in the Circle
\endtitle
\author Julius L. Shaneson \endauthor
\endtopmatter

\centerline{Introduction}
\medskip

Let  $P(t)$ be the number of lattice points in the circle $x^2 + y^2 \le t\,.$
 Gauss showed $P(t) - \pi t = O(t^{1\over 2})\,.$
  Sierpinski decreased the exponent in Gauss's result
to 1/3 ; this followed work of Voronoi on the divisor problem. Van der Corput reduced the exponent to below 1/3. There followed other improvements, Iwaniec and Mozzochi reached 7/22, then further improvements by M.H.Huxley; see the works of this author for more details on the history.

This short paper provides some further estimates of this error, the best of which is
$$P(t) = \pi t + O(t^{1507\over 4875}\ln t)\,\,.$$
By classical methods, a certain sum is transformed up to error into an  exponential sum similar to one also seen classically in the circle problem. This expression is then reformulated, up to error, in a way that permits cancelation of some quadratic terms in an expansion. Some applications of  iterated Weyl-Van der Corput results to the restrictions of the new expression to various ranges then provide the estimate.

For the larger error estimate
$$P(t) = \pi t + O(t^{5\over 16}\ln t)\,\,,$$
the reader can completely omit sections 4 and 6. For the exponent ${509\over 1640}\,,$ section 4 only up to (4.8) is needed and the rest of section 4 and also section 6 can be omitted. For the exponent ${3393\over 10936}\,,$ all of section 4 is needed but not section 6. In all cases, the logarithmic factor can be replaced by a logarithmic factor with a fractional exponent.

In 2007 Cappell and the author posted a paper on the arXiv claiming to obtain the estimate $O(t^{{1\over 4}+ \epsilon})\,.$ Unfortunately we have not been able to produce an error free version. The present paper shares with that paper the Proposition in section 5 and there is also something there akin to what immediately follows the Proposition.
\eject
\centerline{1. A certain sum}
\bigskip
 Let $r$ be a positive integer. A standard result, for example [T,4.7] (compare [GK 3.5]), implies
$$\sum_{m=0}^{[\![\sqrt{t/2}]\!]}{e(2\pi r\sqrt{t-m^2})\over r} =  \sum_{\nu =0}^r{1\over r}\int_{0}^{[\![\sqrt{t/2}]\!]}e(r\sqrt{t-x^2} +\nu x)\,dx+ O(1) \eqno(1.1)$$
(Actually the error could be $O(r^{\epsilon - 1})$ for $\epsilon > 0\,.$)
The critical point of the phase function of the integral on the right for a given $\nu$ is
$$x_\nu = x_{\nu,r,t}= {\nu\sqrt t\over \sqrt{r^2 + \nu^2}}\,.\eqno(1.2)$$
The phase function has the second derivative $-{rt\over y^3}\,,$ where $y = \sqrt{t-x^2}\,.$ Let

$$y_\nu = y_{\nu,r,t}= {r\sqrt t\over \sqrt{r^2 + \nu^2}} = \sqrt{t-x_\nu^2}\,.\eqno(1.3)$$

The van der Corput estimate (e.g. [GK 3.2][T,4.4]) gives

$$\int_{0}^{[\![\sqrt{t/2}]\!]}e(r\sqrt{t-x^2} -\nu x) dx \ll { t^{1\over 4}\over r^{1\over 2}}\,\,\,,\,\,\,0 \le \nu \le r\eqno (1.4)$$
We have ${\partial x_\nu\over \partial x_\nu} = {r^2\sqrt t\over (r^2+\nu^2)^{3\over 2}}\ge {\sqrt t\over 2\sqrt 2 r}\,,$ $0\le \nu \le r\,;$ therefore for $1 \le \nu \le r-1$
$${\tau\sqrt t\over 2\sqrt 2 r} \le x_\nu \le \sqrt {t\over 2} - {\tau\sqrt t\over 2\sqrt 2 r}\,.\eqno(1.5)$$
Assuming further that $r \le {t^{1\over 2}\over 2} \,,$ it then follows from [GK,3.4] (compare [T,4.6]) that for
$1 \le \nu \le r-1\,,$

$$\int_{0}^{[\![\sqrt{t/2}]\!]}e(r\sqrt{t-x^2} +\nu x) = {e(-1/8)rt^{1\over 4}\over (r^2 + \nu^2)^{3\over 4}}e(ry_\nu + \nu x_\nu)+O(1)\,.\eqno(1.6)$$
Therefore

$$\sum_{m=0}^{[\![\sqrt{t/2}]\!]}{e(2\pi r\sqrt{t-m^2})\over r} = \sum_{\nu = 0}^r{e(-1/8)t^{1\over 4}\over (r^2 + \nu^2)^{3\over 4}}e(ry_\nu + \nu x_\nu) + O\Big( { t^{1\over 4}\over r^{3\over 2}} + 1)\,.\eqno(1.7)$$
\eject
\centerline{2. A reformulation}
\bigskip

For $1 \le \nu \le r\,,$ let $m_\nu \in Z[{1\over \nu}]$ the be largest rational number with denominator $\nu$ that is at most $x_{\nu}\,.$ Let $m_0 = 0\,.$As above, write $y= y(x) = \sqrt{t-x^2}\,.$  Let $n_\nu = \sqrt{t-m_\nu^2}\,.$
Let $\lambda_\nu = x_\nu - m_\nu\,.$
Then $y' = -{x\over y}$ and $y'' = -ty^{-3}\,.$ So $n_0 = y_0 = \sqrt t$ and for $1 \le \nu \le r\,,$

$$ n_\nu = y_\nu + {\nu\over r}\,\lambda_\nu  + O(t^{-{1\over 2}}\nu^{-2})\,,\eqno(2.1)$$
from which it follows that

$$e(ry_\nu + \nu x_\nu) = e(ry_\nu + \nu\lambda_\nu) = e(r(y_\nu + {\nu\over r}\,\lambda_\nu))= e(rn_\nu) + O\Big({r\over \nu^2\sqrt t}\Big)\,,\eqno(2.2)$$
So, assuming $r \le t^{1\over 2}\,.$

$$\sum_{m=0}^{[\![\sqrt{t/2}]\!]}{e(2\pi r\sqrt{t-m^2})\over r} = \sum_{\nu = 0}^r{e(1/8)t^{1\over 4}\over (r^2 + \nu^2)^{3\over 4}}e(rn_\nu) + O\Big({t^{1\over 4}\over r^{3\over 2}} + 1)\,;\eqno(2.3)$$
the error in (2.2) gets absorbed in that of (1.7).

For a fixed $1 < c < \sqrt 2\,,$  there exist positive constants $c_i < C_i$ (that depend on $c$ but not $t\,$)
so that $c_it^{-{i-1\over 2}} \le |y^{(i)}(x)| \le C_it^{-{i-1\over 2}}\,,$ for $|x| \le c^{-1}\sqrt t\,;$ i.e. $y^{(i)}(x) \approx t^{-{i-1\over 2}}$ in this range.
It is also not difficult to see that for $0 \le \nu \le r-1\,,$

$${\sqrt t\over 2\sqrt 2 r} \le x_{\nu + 1} - x_\nu \le {\sqrt t\over r}\,,\eqno(2.4)$$
the extremes being lower and upper bounds of ${\partial x_\nu\over \partial\nu}\,,$ $0 \le \nu \le r\,.$
Expanding around $x_{\nu}\,,$

$$n_{\nu}-n_{\nu-1} = -{\nu\over r}(m_{\nu } - m_{\nu-1}) + \sum_{i=2}^{N-1}{y^{(i)}(x_{\nu})\over i\,!}\Bigg[(-\lambda_\nu)^i-(m_{\nu-1} - x_\nu)^i\Bigg] + O\Big({\sqrt t\over r^N}\Big)\,\eqno(2.5.1)$$
for $1 \le \nu \le r\,,$ and for $0 \le \nu \le r-1$

$$\eqalign{ n_{\nu+1}-n_{\nu} &= \!-{\nu\over r}(m_{\nu +1} - m_{\nu}) + \! \sum_{i=2}^{N-1}\!{y^{(i)}(x_{\nu})\over i\,!}\Bigg[(m_{\nu+1} - x_{\nu})^i  -(-\lambda_{\nu})^i\Bigg] \cr&+ O\Big({\sqrt t\over r^N}\Big)\,.\cr}\eqno(2.5.2)$$

The sum of these will give an expression for the difference $n_{\nu+1}-n_{\nu-1}\,;$ i.e. up to the errors this difference will be the sum of a rational number with denominator $r$ and

$$G(\nu) = \sum_{i=2}^{N-1}{y^{(i)}(x_\nu)\over i\,!}\Bigg[(m_{\nu+1} - x_\nu)^i -(m_{\nu-1} - x_{\nu})^i\Bigg]\,.\eqno(2.6)$$
Let

$$F(\nu) = n_1 + \sum_{\mu = 1}^{\nu-1\over 2}G(2\mu), \,\,\,\nu\,\,\, \hbox{odd}\,,\,\,\,F(\nu) = n_0 + \sum_{\mu = 1}^{\nu\over 2}G(2\mu - 1),\,\,\,\nu\,\,\, \hbox{even}\,.\eqno(2.7)$$
($F(0) = n_0\,,F(1) = n_1)\,.$)
Then from (2.5.1) and (2.5.2)

$$e(rn_\nu) = e(rF(\nu)) + O\Big( {\nu\sqrt t\over r^{N-1}}\Big)\,.\eqno(2.8)$$
Hence (1.7) becomes

$$\sum_{m=0}^{[\![\sqrt{t/2}]\!]}{e(2\pi r\sqrt{t-m^2})\over r} = \sum_{\nu = 0}^r{e(-1/8)t^{1\over 4}\over (r^2 + \nu^2)^{3\over 4}}e(rF(\nu)) + O\Big({t^{1\over 4}\over r^{3\over 2}} + 1 + { t^{3\over 4}\over r^{N-{3\over 2}}} \Big)\,.\eqno(2.9)$$

Let

$$G_1(\nu) = \sum_{i=2}^{N-1}{y^{(i)}(x_\nu)\over i\,!}\Bigg[(x_{\nu+1} - x_\nu)^i -(x_{\nu-1} - x_{\nu})^i\Bigg]\,.\eqno(2.10)$$
Then it is not difficult to see, using (2.4) and $0 \le \lambda_\nu \le \nu^{-1}\,,$  that
$$G(\nu) - G_1(\nu) = O((r\nu)^{-1})\,.\eqno(2.11)$$
\big(More definitely, one can see that there is a constant $C$ (depending on $N$) such that

$$ 0 < G(\nu) \le G_1(\nu) + {C\over r\nu}\,.\Bigg)$$
Therefore, let $F_1(\nu)$ be defined as $F$ is in (2.7), but using $G_1$ in place of $G\,.$
Then there is an expression

$$\sum_{m=0}^{[\![\sqrt{t/2}]\!]}{e(2\pi r\sqrt{t-m^2})\over r} = \sum_{\nu = 0}^r{e(1/8)t^{1\over 4}\over (r^2 + \nu^2)^{3\over 4}}\,e(\omega_\nu)e(rF_1(\nu)) + O\Big({t^{1\over 4}\over r^{3\over 2}} + 1 + { t^{3\over 4}\over r^{N-{3\over 2}}} \Big)\,,\eqno(2.12)$$
in which the coefficients $e(\omega_\nu)$ have total variation $O(\ln r)\,;$
that is,

$$\sum_{\nu = 1}^r|\,e(\omega_\nu) - e(\omega_{\nu-1})| \le \sum_{\nu= 1}^r {C\over \nu} \ll \ln  r\,.\eqno(2.13)$$

Finally, it is not hard to see there is a positive constants  $D_i$ such that
for $0 \le \nu \le r\,,$

$$ \Big| {\partial^ix_{\nu}\over \partial \nu^i}\Big| \le {D_i\sqrt t\over r^i}\,.\eqno(2.14)$$
Therefore, for $\delta = \pm 1\,$ and $1 \le \nu \le r-1\,,$

$$x_{\nu+\delta} - x_\nu = \sum_{k=1}^{M-1}{1\over k!}{\partial^kx_{\nu}\over \partial \nu^k}\,\delta^k + O\Big({\sqrt t\over r^M}\Big)\,.\eqno(2.15)$$
It then follows that $G_1$ can replaced by

$$G_2(\nu) = \sum_{i=3\atop i odd}^{N-1}{2y^{(i)}(x_\nu)\over i\,!}\Bigg[\sum_{k=1}^{N-i}{1\over k!}{\partial^kx_{\nu}\over \partial \nu^k}\Bigg]^i\eqno(2.16)$$
without further error. That is, if $F_2(\nu)$ is defined as in (2.7) with $G_2$ instead of $G$ or $G_1\,,$ then (2.12) still holds, i.e.

$$\sum_{m=0}^{[\![\sqrt{t/2}]\!]}{e(2\pi r\sqrt{t-m^2})\over r} = \sum_{\nu = 0}^r{e(-1/8)t^{1\over 4}\over (r^2 + \nu^2)^{3\over 4}}\,e(\omega_\nu)e(rF_2(\nu)) + O\Big({t^{1\over 4}\over r^{3\over 2}} + 1 + { t^{3\over 4}\over r^{N-{3\over 2}}} \Big)\,.\eqno(2.17)$$
\bigskip

\centerline{3.  An Estimate}
\bigskip
From the previous section, $y^{(i)}(x) \approx -t^{-{i-1\over 2}}$ for $i \ge 2\,,$  $|x| \le c^{-1}\sqrt t$ and
${\partial^ix_{\nu}\over \partial \nu^i} \ll \sqrt t r^{-i}\,.$
Therefore, it is not hard to see,

$${\partial^l G_2(\nu)\over \partial \nu^l} = {1\over 3} {\partial^l\over \partial\nu^l}\Bigg[y^{(3)}(x_\nu)\Bigg({\partial x_\nu\over \partial\nu}\Bigg)^3\,\Bigg] + A_l\,,\eqno(3.1)$$
with
$$A_l \ll {\sqrt t\over r^{l+4}}\,.\eqno(3.2)$$
More explicitly,

$${\partial^l G_2(\nu)\over \partial \nu^l} = - {\partial^l\over \partial\nu^l}\Bigg[{r\nu\sqrt t\over (r^2 + \nu^2)^{5\over 2}}\Bigg] + A_l\,,\eqno(3.3)$$
In particular

$${\partial^2 G_2(\nu)\over \partial \nu^2} = -{5r\nu(4\nu^2 - 3r^2)\sqrt t\over (r^2 + \nu^2)^{9\over 2}} + A_2\,.\eqno(3.4)$$

Consider the sums, $\nu$ ranging from 1 to some $r_1 \le r\,,$

$$S_1 = \sum_{\nu\,odd}e(rF_2(\nu))\qquad\qquad S_2 = \sum_{\nu\,even}e(rF_2(\nu))\,.\eqno(3.5)$$
We now want to apply Theorem 2.6 of [GK] to estimate these sums. To do so, let $t^{\alpha} < \omega < (1 - {\sqrt 3\over 2})r$ for a fixed $\alpha < 1\,.$ Let $\nu_0 = {\sqrt 3\over 2}r\,.$ Let
$I_0 = [1,\omega]\,,$ $I_1 = [\omega, \nu_0 - \omega)\,,$ $I_2 = [\nu_0 - \omega,\nu_0 + \omega]\,,$ $I_3 = (\nu_0 + \omega, r]\,.$ (The upper bound on $\omega$ guarantees that this is a increasing sequence of non-empty disjoint intervals.) On the intervals $I_1$ and $I_3\,,$ from calculus methods

$$ {5(\sqrt 3 -1)\over 4\sqrt t}{\omega\sqrt t\over r^6} \le \Big|{r\nu(4\nu^2 - 3r^2)\over (r^2 + \nu^2)^{9\over 2}}\Big| \le 5 {\sqrt t\over r^5} = 5{\omega\sqrt t\over r^6}\cdot{r\over \omega}\,.\eqno(3.6)$$
(The actual constants not important and these are not the best.)
Therefore, from (3.2) and the lower bound on $\omega\,,$ on these intervals

$${\omega\sqrt t\over r^6} \ll \Big|{\partial^2 G_2(\nu)\over \partial \nu^2}\Big|
\ll {\omega\sqrt t\over r^6}\cdot{r\over \omega}\,,\eqno(3.7)$$
with the implied constants depending only on $\alpha$ and the constant in (3.2).
Let, $\nu$ still confined to the range between 1 and $r_1\,,$

$$ \tilde S_1 = \sum_{\nu\,odd\atop \nu\in I_1\cup I_3}e(rF_2(\nu))\qquad\qquad \tilde S_2 = \sum_{\nu\,even\atop \nu\in I_1\cup I_3}e(rF_2(\nu))\,.\eqno(3.8)$$
Then 2.6 of [GK], with the second derivative of $G_2$ in the role of the third derivative of the phase function, would state:

$$\tilde S_i \ll {t^{1\over 12}r^{1\over 2}\over \omega^{1\over 6}} + {r\over \omega^{1\over 4}} + {r^{3\over 2}\over t^{1\over 8}\omega^{1\over 4}}\,.\eqno(3.9)$$

This estimate would follow precisely from (3.7) and (2.6) of [GK] if $G_2$ were the derivative of
of $F_2\,.$ However, there is just the relationship $G_2(\nu) = F_2(\nu+1) - F_2(\nu-1)\,.$
It seems clear that  the proof of Theorem 2.6 of [GK] goes through almost unchanged in this context. The only significant difference would be that Theorem 2.1 in [GK] as used in the proof there needs to be replaced  by the stronger Kusmin-Landau inequality, which is stated for example in section 2.6 of [GK].

To be able to apply the results of [GK] exactly as stated and proved, a suitable version of the Euler-MacLaurin formula can be used as follows.  Let
$I$ be an interval with first least odd integer $\nu_1$ in $I_1$ or $I_3$, and consider a sum, for example,

$$ S = \sum_{\nu\,odd\atop \nu\in I}e(rF_2(\nu))\,.\eqno(3.10)$$
Let

$$ F_3(\nu) = {1\over 2}\int_{\nu_1}^{\nu}G_2(u)du + \sum_{j=2} ^N\bar b_{j}\Big(G^{(j-1)}(\nu) - G^{(j-1)}(\nu_1)\Big)\eqno(3.11)$$
Here $\bar b_{j} = {2^j\bar B_j\over j!}\,,$ $\bar B_j$ the values of  Bernoulli polynomials at ${1\over 2}$ (actually zero for $j$ odd). From (3.3), ${\partial^N G(\nu)\over \partial\nu^N} \ll {\sqrt t\over r^{N+3}}\,.$ Therefore the Euler-MacLaurin formula with remainder implies that for $\nu$ odd

$$F_2(\nu) = F_2(\nu_1) + F_3(\nu) + O\Big({\sqrt t\over r^{N +2}}\Big)\,.\eqno(3.12)$$

Hence

$$S = e(rF_2(\nu_1))\sum_{\nu\,odd\atop \nu\in I}e(rF_3(\nu)) + O\Big({\sqrt t\over r^{N}}\Big)\,.\eqno(3.13)$$
The derivative of order $l+1$ of $F_3(\nu)$ is one-half the derivative of $G_2(\nu)$ of  order $l$ plus a linear combination of higher order ones. The higher order derivatives can be absorbed in $A_l; $ thus (3.3) holds with the derivative of order $l$ of $G_2$ replaced by derivative of order $l+1$ of $F_3\,.$ In particular (3.7) also holds for
the third derivative of $F_3$ on $I \subset I_1\cup I_3\,.$  Therefore, $S$  satisfies the conclusion of Theorem (2.6) of [GK], at least possibly up to the further error (a larger error has already appeared) of $O(\sqrt t r^{-N})\,.$ The same type of argument works for sums over even integers.

Combining (3.9) with the trivial estimate on the intervals $I_0$ and $I_2$  gives

$$S_i \ll \omega + {t^{1\over 12}r^{1\over 2}\over \omega^{1\over 6}} + {r\over \omega^{1\over 4}} + {r^{3\over 2}\over t^{1\over 8}\omega^{1\over 4}}\,.\eqno(3.14)$$
Note that the third term is at most the second if $r \le t^{1\over 4}\,.$
Set $\omega = t^{1\over 14}r^{3\over 7}\,.$ Then, as long as $(1 -{\sqrt 3\over 2})^{-{7\over 4}} t^{{1\over 8}} < r \le t^{1\over 4}\,,$

$$S_i \ll t^{1\over 14}r^{3\over 7}  + {r^{25\over 28}\over t^{1\over 56}}  \,.\eqno(3.15)$$
This estimate can be applied in (2.17), using (2.13) and Abel summation. Taking a large enough enough $N\,,$ it results that for the same range of $r\,,$

$$\sum_{m=0}^{[\![\sqrt{t/2}]\!]}{e(2\pi r\sqrt{t-m^2})\over r} \ll
{t^{9\over 28}\ln r\over r^{15\over 14}} + {t^{13\over 56}\ln r\over r^{17\over 28}} + {t^{1\over 4}\over r^{3\over 2}} + 1\,.\eqno(3.16)$$
Hence, for $(1 -{\sqrt 3\over 2})^{-{7\over 4}}t^{{1\over 8}} < R_1 < R \le t^{1\over 4}\,,$

$$\sum_{r=R_1}^R \sum_{m=0}^{[\![\sqrt{t/2}]\!]}{e(2\pi r\sqrt{t-m^2})\over r} \ll
{t^{9\over 28}\ln R\over R_1^{1\over 14}} + t^{13\over 56}R^{11\over 28}\ln R   + {t^{1\over 4}\over R_1^{1\over 2}} + R\,.\eqno(3.17)$$

.
On the other hand, the trivial estimate on the sum on the right in (1.8) gives (eliminating superfluous terms)

$$\sum_{r=1}^{R_1} \sum_{m=0}^{[\![\sqrt{t/2}]\!]}{e(2\pi r\sqrt{t-m^2})\over r} \ll
t^{1\over 4}R_1^{1\over 2} + R_1\,.\eqno(3.18)$$
So, taking for $R_1$ the least integer greater than  $(1 -{\sqrt 3\over 2})^{-{7\over 4}}t^{{1\over 8}}\,$ gives

$$\sum_{r=1}^{R} \sum_{m=0}^{[\![\sqrt{t/2}]\!]}{e(2\pi r\sqrt{t-m^2})\over r} \ll
t^{{5\over 16}}\ln t + t^{13\over 56}R^{11\over 28}\ln R  \,,\eqno(3.19)$$
valid for $1 \le R \le t^{1\over 4}\,.$

\bigskip
\centerline{4. Further estimates}
\bigskip

To obtain an improvement, we apply to the sums (3.5) the next case $q=2$ of Theorem (2.8) of [GK].
There is  explicit formula

$${\partial^3 G_2(\nu)\over \partial \nu^3} = {15r\sqrt t(8\nu^4 - 12\nu^2r^2 + r^4)\over (r^2 + \nu^2)^{11\over 2}} + A_3\,.\eqno(4.1)$$
The only value in $ [0,r]$ for which this derivative vanishes is $\eta = {(3-\sqrt 7)^{1\over 2}\over 2}\,r = (1 - {\sqrt 7\over 3})^{1\over 2}\,\nu_0\,.$

Let  $J_1 = [0,\eta - \rho),,$ $J_2 = [\eta-\rho,\eta + \rho]\,,$ and $J_3 = (\eta + \rho,r]\,.$ Then, similarly (slightly harder) to (3.7),

$${\rho\sqrt t\over r^7} \ll \Big|{\partial^3 G_2(\nu)\over \partial \nu^3}\Big|
\ll {\rho\sqrt t\over r^7}\cdot{r\over \rho}\eqno(4.2)$$
on the intervals $J_1$ and $J_3\,,$ assuming  $t^\alpha < \rho < {(3-\sqrt 7)^{1\over 2}\over 2}\,r\,.$
Let

$$ \bar S_1 = \sum_{\nu\,odd\atop \nu\in J_1\cup J_3}e(rF_2(\nu))\qquad\qquad \bar S_2 = \sum_{\nu\,even\atop \nu\in J_1\cup J_3}e(rF_2(\nu))\,.\eqno(4.3)$$
Then the case $q =2$ an $Q =4$ of Theorem (2.8) of [GK] asserts:

$$\bar S_i \ll {t^{1\over 28}r^{5\over 7}\over \rho^{1\over 14}} + {r\over \rho^{1\over 8}} + {r^{21\over 16}\over \rho^{1\over 8}t^{1\over 16}}\,.\eqno(4.4)$$

Combining this with the trivial estimate over $J_2$ gives

$$ S_i \ll \rho + {t^{1\over 28}r^{5\over 7}\over \rho^{1\over 14}} + {r\over \rho^{1\over 8}} + {r^{21\over 16}\over \rho^{1\over 8}t^{1\over 16}}\,.\eqno(4.5)$$
Note that the last term is no larger than the previous one for $r \le t^{1\over 5}\,.$
Set $\rho = t^{1\over 30}r^{2\over 3}\,.$ Then, as long as ${8\over (3-\sqrt 7)^{3\over 2}}t^{1\over 10} < r \le t^{1\over 5}\,,$

$$S_i \ll t^{1\over 30}r^{2\over 3} + {r^{11\over 12}\over t^{1\over 240}} \eqno(4.6)$$
Hence, for $ {8\over (3-\sqrt 7)^{3\over 2}}t^{{1\over 10} } < R_2 < R_1 < t^{1\over 5}\,,$

$$\sum_{r=R_2}^{R_1} \sum_{m=0}^{[\![\sqrt{t/2}]\!]}{e(2\pi r\sqrt{t-m^2})\over r} \ll t^{17\over 60}R_1^{1\over 6}\ln R_1 + t^{59\over 240}R_1^{5\over 12}\ln R_1 + {t^{1\over 4}\over R_2^{1\over 2}} + R_1\,.\eqno(4.7)$$

This can be combined with (3.17) and with the trivial estimate for $1 \le r \le R_2\,.$ We take $R_1$ to be the nearest integer to the minimum of $R$ and $t^{127\over 820}\,;$ this value is obtained by equating the second term in (4.7) with the first term in (3.17), neglecting the logarithmic terms. For $R_2$ we take the least integer more than ${8\over (3-\sqrt 7)^{3\over 2}}t^{1\over 10}\,.$
After eliminating superfluous terms, the result is

$$\sum_{r=1}^{R} \sum_{m=0}^{[\![\sqrt{t/2}]\!]}{e(2\pi r\sqrt{t-m^2})\over r} \ll
t^{{509\over 1640}}\ln t + t^{13\over 56}R^{11\over 28}\ln R \,.\eqno(4.8)$$
for $R \le t^{1\over 5}\,.$
\bigskip

For a further improvement, we try to use a non trivial estimate on
$$ \check S_1 = \sum_{\nu\,odd\atop \nu\in J_2}e(rF_2(\nu))\qquad\qquad \check S_2 = \sum_{\nu\,even\atop \nu\in J_2}e(rF_2(\nu))\,.\eqno(4.9)$$
Explicitly,

$${\partial^4 G_2(\nu)\over \partial \nu^4} = -{105r\nu\sqrt t(8\nu^4 -20\nu^2r^2 + 5r^4)\over (r^2 + \nu^2)^{13\over 2}} + A_4\,.\eqno(4.10)$$
The first term on the right is zero in the relevant range only at

$$\nu = {(5 -\sqrt 15)^{1\over 2}\over 2}\, r\,,$$
which is distinct multiple of $r$ from $\eta = {(3-\sqrt 7)^{1\over 2}\over 2}\,r\,.$ Another calculus type argument shows that, as long as  $\rho < kr$ ($k$ say half the difference of the two above multipliers of $r$),

$${\partial^4 G_2(\nu)\over \partial \nu^4} \approx {\sqrt t\over r^7}\eqno(4.11)$$
on the interval $J_2\,.$
Then Theorem (2.8) of [GK] with $q =3$ and $Q = 8$ asserts

$$\check S_i \ll {\rho t^{1\over 60}\over r^{1\over 5}} + \rho^{15\over 16} + {\rho^{49\over 64}r^{3\over 8}\over t^{1\over 32}}\,.\eqno(4.12)$$
Therefore, if $t^\alpha < \rho < kr\,,$ and $r \le t^{1\over 5}\,,$

$$S_i = \bar S_i + \check S_i \ll {t^{1\over 28}r^{5\over 7}\over \rho^{1\over 14}} + {r\over \rho^{1\over 8}}  + {\rho t^{1\over 60}\over r^{1\over 5}} + \rho^{15\over 16} + {\rho^{49\over 64}r^{3\over 8}\over t^{1\over 32}}\,.\eqno(4.13)$$

Set $\rho = t^{2\over 25}r^{152\over 375}\,;$ this equalizes the first and last terms. Then, as long as $ k_1t^{30\over 223} < r \le  t^{1\over 5}\,,$ for a suitable constant $k_1,$

$$S_i \ll t^{3\over 100}r^{257\over 375} +{r^{356\over 375}\over t^{{1\over 100}}} + t^{29\over 300} r^{77\over 375} + t^{3\over 40}r^{19\over 50}\,.\eqno(4.14)$$
Hence, for $0  < k_1t^{{30\over 223} } < R_2 < R_1 < t^{1\over 5}\,,$

$$\eqalign{\sum_{r=R_2}^{R_1} \sum_{m=0}^{[\![\sqrt{t/2}]\!]}{e(2\pi r\sqrt{t-m^2})\over r} &\ll t^{7\over 25}R_1^{139\over 750}\ln R_1 + t^{6\over 25}R_1^{337\over 750}\ln R_1+ {t^{52\over 150}\ln R_1\over R_2^{221\over 750}}  \cr&+
{t^{13\over 40}\ln R_1\over R_2^{3\over 25}}  + {t^{1\over 4}\over R_2^{1\over 2}} + R_1\,.}\eqno(4.15)$$

This can be combined with (3.17). Take for $R_1$ the nearest integer to $t^{855\over 5648}$ ; this results from equating the second term in (4.15) and the first of (3.17), neglecting logarithmic terms. With the substitution, the first term in (4.15) is slightly smaller than the second (equal to the last), and the last two terms  of (3.17) and (4.15) are also superfluous, so for
$0  < k_1t^{{30\over 223} } < R_2 < R < t^{1\over 5}\,,$

$$\sum_{r=R_2}^{R} \sum_{m=0}^{[\![\sqrt{t/2}]\!]}{e(2\pi r\sqrt{t-m^2})\over r} \ll
t^{3393\over 10936}\ln t + t^{13\over 56}R^{11\over 28}\ln R + {t^{52\over 150}\ln R_1\over R_2^{221\over 750}} + {t^{13\over 40}\ln R_1\over R_2^{3\over 25}}\,.\eqno(4.16)$$

Now take $R_2$ to be the first integer more than $k_1t^{30\over 223}\,,$ apply (4.7) in the range $t^{1\over 10} < r \le k_1t^{30\over 223}$ and the trivial estimate below this range. All the terms obtained, i.e. the last two in (4.16) and what we get from (4.7) and the trivial estimate, are majorized by the first term of (4.16). Therefore

$$\sum_{r=1}^{R} \sum_{m=0}^{[\![\sqrt{t/2}]\!]}{e(2\pi r\sqrt{t-m^2})\over r} \ll
t^{3393\over 10936}\ln t + t^{13\over 56}R^{11\over 28}\ln R  \,,\eqno(4.17)$$
$1 \le R \le t^{1\over 5}\,.$

\bigskip
\centerline{5. Estimates on $P(t) - \pi t\,.$}
\bigskip

Let $\psi(t) = t - [\![t]\!] - {1\over 2}$ be the usual saw-tooth function.
This section first (re)states analogues of  well known results for the divisor problem.
\proclaim{ Proposition }

$$ \sum_0^{[\![\sqrt{t/2}]\!]} \psi(\sqrt{t - m^2}) = {\pi t - P(t)\over 8} + O(1)\,\,.\eqno(5.1)$$

\endproclaim
\noindent
Proof:\smallskip

 By trapezoidal approximation using integer points from $0$ to the greatest integer $\alpha$ in $\sqrt{t\over 2}\,\,,$

$$\sum_0^{\alpha} \psi(\sqrt{t - m^2}) = A - L + {\sqrt t\over 2} + O(1)\,\,,\eqno(5.2)$$
where $A$ is the area under the circle of radius $y  = \sqrt{t-x^2} $ for $0 \le x \le \alpha$
and $L$ the number of lattice points in this region, counting the top and sides but not the bottom. Combining this with "Pick's theorem"

$$L_0 = A_0 + {3\over 2}\alpha + 1\eqno(5.3)$$
for the number of lattice points $L_0$ in the triangle with vertices the origin, $(\alpha,0)\,\,,$ and $(\alpha,\alpha)\,\,,$ and with area $A_0\,,$ it follows that

$$\sum_1^{\alpha} \psi(\sqrt{t - m^2}) = {\pi t\over 8} - L_1 + {\sqrt t\over 2} - {1\over 2}\sqrt{t\over 2} + O(1)\,\,,\eqno(5.4)$$
$L_1$ the number of lattice points contained in the sector of the circle of radius
$\sqrt t$ with angle ${\pi\over 4} < \theta \le {\pi\over 2}\,\,.$
However, it is obvious that the number of lattice points in the circle is

$$P(t) = 8L_1 + 4\sqrt{t\over 2} - 4\sqrt t\,,\eqno(5.5)$$
and the result follows.
\bigskip

Let $\Vert u \Vert$ be the distance from a real number $u$ to the nearest integer.
The Fourier expansion of $\psi$ gives

$$\sum_0^{\sqrt{t/2}} \psi(\sqrt{t - m^2}) = -{1\over \pi}\sum_0^{\sqrt{t/2}}\sum_{r=1}^\infty{\sin 2\pi r\sqrt{t-m^2}\over r}\,,\eqno(5.3)$$
Let $J $ denote the set of lattice points $(m,n)$ with $0 \le m \le \alpha\,,$ ($\alpha$ as in the previous proof) and  with $\sqrt{t-m^2} - {1\over 2} \le n < \sqrt{t-m^2} + {1\over 2}\,.$
Then

$$\eqalign{\sum_0^{\sqrt{t/2}}\sum_{r=R+1}^\infty{\sin 2\pi r\sqrt{t-m^2}\over r} &\ll \sum_0^{\sqrt{t/2}} \min\Big\{{\Vert\sqrt{t-m^2}\Vert^{-1}\over R}, 1\Big\}\cr&\le \sum_{(m,n)\in J} \min\Big\{{1\over R|\sqrt{m^2 + n^2}-\sqrt t\,|}\,, 1\Big\} \,,\cr}\eqno(5.4)$$
The last inequality follows from the fact the closest point on a circle to a given point is the intersection with the circle of the line joining that point to the center of the circle. This fact also means that if $J_k$ consists of those points in $J$ with
${k\over 2R} \le \sqrt{m^2+n^2}-\sqrt t < {k+1\over 2R}\,,$ then
$$J = \bigcup_{k = -R}^{R-1} J_k\,.\eqno(5.5)$$
Therefore

$$\sum_0^{\sqrt{t/2}}\sum_{r=R+1}^\infty{\sin 2\pi r\sqrt{t-m^2}\over r}\ll
\sum_{k=-R}^{-2}{ n(J_k)\over -k-1}+ \sum_{k=2}^{R-1}{ n(J_k)\over k} + n(J_{-1}) + n(J_0)\,,\eqno(5.6)$$
$n(J_k)$ the number of elements of $J_k\,.$
Further, distinct circles of radius less than $\sqrt s$ that contain
lattice points must differ in radius by at least ${1\over 2\sqrt s}\,;$
as is well known, each such circle has at most $O(s^\epsilon)$ lattice points, any $\epsilon > 0\,.$ Therefore $n(J_k) \ll t^{{1\over 2} + \epsilon}R^{-1}\,,$ assuming $R \le t^{1\over 2}\,, $ and so, also using $\ln t \ll t^{\epsilon}\,,$

$$\pi t - P(t) = -{8\over \pi}\sum_{m=0}^{\sqrt{t\over 2}}\,\sum_{r=1}^{R}{\sin 2\pi r\sqrt{t-m^2}\over r}
+ O({t^{{1\over 2} + \epsilon}\over R})\,.\eqno(
5.7)$$

Therefore, for $R < t^{1\over 5}\,,$

$$P(t) - \pi t \ll t^{\beta}\ln t + t^{13\over 56}R^{11\over 28}\ln R +  {t^{{1\over 2} + \epsilon}\over R}\,.\eqno(5.8)$$
with $\beta = {5\over 16}$ from (3.19), $\beta = {509\over 1640}$ from (4.8), or
$\beta = {3393\over 10936}$ from (4.17)
Choose $\epsilon$ small and $R$ to be the nearest integer to $t^{15\over 78}\,;$ then all the other terms are less than the first, i.e.

$$P(t) - \pi t \ll t^{\beta}\ln t\,.\eqno(5.9)$$

\medskip
\centerline{6. A further estimate}
\medskip

In place of the trivial estimate on the intervals $I_0$ and $I_2$ used to obtain (3.14), observe that for $  \omega < k_2r\,,$ for a suitable constant $k_2\,,$

$${\partial^3 G_2(\nu)\over \partial \nu^3} \approx {\sqrt t\over r^6} \eqno(6.1)$$
on these intervals. Therefore, we can apply Theorem (2.8) of [GK] with $q = 2$ to replace (3.14) with

$$S_i \ll {\omega t^{1\over 28}\over r^{5\over 14}} + \omega^{7\over 8} + {\omega^{9\over 16}r^{5\over 8}\over t^{1\over 16}} + {t^{1\over 12}r^{1\over 2}\over \omega^{1\over 6}} + {r\over \omega^{1\over 4}} + {r^{3\over 2}\over t^{1\over 8}\omega^{1\over 4}}\,,\eqno(6.2)$$
$t^\alpha < \omega < k_3 r\,.$
Set $\omega = t^{2\over 25}r^{12\over 25}\,,$ the result of equating the second and fourth terms. Then, as long as $ k_4t^{2\over 13} < r \le t^{1\over 4}\,,$ for a suitable $k_4\,,$

$$S_i \ll {t^{81\over 700} r^{43\over 350}}+  t^{7\over 100}r^{21\over 50} + {r^{179\over 200}\over t^{7\over 400}} + {r^{22\over 25}\over t^{1\over 50}} \,.\eqno(6.3) $$
Note the last term is majorized by the previous one and the first by the second.

Therefore, analogous to previous sections, for $k_4t^{2\over 13} < R_1 < R < t^{1\over 4}\,,$

$$\sum_{r=R_1}^R \sum_{m=0}^{[\![\sqrt{t/2}]\!]}{e(2\pi r\sqrt{t-m^2})\over r} \ll {t^{8\over 25}\ln t\over R_1^{2\over 25}}+ t^{93\over 400}R^{79\over 200}\ln t  + {t^{1\over 4}\over R_1^{1\over 2}} + R\,.\eqno(6.4)$$

Now suppose $k_4t^{2\over 13} < R_1 < R < t^{1\over 5}\,,$ $k_2t^{30\over 223} < R_2 < R_1\,.$
Then by (6.4), (4.15), (4.7) in the range ${8\over (3-\sqrt 7)^{3\over 2}}\,t^{1\over 10} < r \le R_2\,,$ and the trivial estimate below this range,

$$\eqalign{\sum_{r=1}^R \sum_{m=0}^{[\![\sqrt{t/2}]\!]}&{e(2\pi r\sqrt{t-m^2})\over r} \ll {t^{8\over 25}\ln t\over R_1^{2\over 25}}+ t^{93\over 400}R^{79\over 200}\ln t + t^{7\over 25}R_1^{139\over 750}\ln R_1 + t^{6\over 25}R_1^{337\over 750}\ln R_1 \cr&+ {t^{52\over 150}\ln R_1\over R_2^{221\over 750}}  +
{t^{13\over 40}\ln R_1\over R_2^{3\over 25}} + t^{17\over 60}R_2^{1\over 6}\ln R_2 + t^{59\over 240}R_2^{5\over 12}\ln R_2 + t^{3\over 10}\,.}$$
Let $R_1$ be the first integer at least $k_4t^{2\over 13}$ and $R_2$ the nearest integer to $t^{95\over 692}\,$ (gotten by equating the fifth and seventh terms). Then
the fourth term majorizes everything else except possibly the second, and we get:

$$\sum_{r=1}^R \sum_{m=0}^{[\![\sqrt{t/2}]\!]}{e(2\pi r\sqrt{t-m^2})\over r} \ll
 t^{1815\over 5876}\ln t+ t^{93\over 400}R^{79\over 200}\ln t \,. \eqno(6.5)$$
Note that if using (6.4) with a smaller lower bound could be justified, the estimate could be improved further.
In any case, as in the previous section,

$$P(t) - \pi t \ll t^{1507\over 4875}\ln t + t^{93\over 400}R^{79\over 200}\ln t  + {t^{{1\over 2} + \epsilon}\over R}\,.\eqno(6.6)$$
Take $R$ the nearest integer to $t^{107\over 558}$ and $\epsilon$ small. Then  the other terms are less than the first, i.e.

$$P(t) - \pi t \ll t^{1507\over 4875}\ln t\,.\eqno(6.7)$$
\bigskip
\centerline{Concluding remark.}
\bigskip

The method of exponent pairs, as presented in [GK] and elsewhere, does not apply directly to the exponential sums considered here. Rather, results like those quoted produce exponent systems as originally defined by Van der Corput. However, methods used to get various exponents pairs should able to be extended to produce further exponent systems beyond the ones quoted here and thereby possibly also produce better estimates for the circle problem.

\bigskip

\medskip
\centerline{References}
\medskip

\smallskip
\noindent
[Gaus1] C.F. Gauss, {\it De nexu inter multitudinem classium in quas formae
binariae secundi gradus distribuuntur earumque determinantem.} Werke (1863),
Vol 2, pp. 269-291 (especially, p. 280)
\smallskip
\noindent
[Gaus2] C.F. Gauss, {\it Disquisitiones arithmeticae.} (German edition by H. Maser). (1886) p. 657
\smallskip
\noindent
[GK] S.W. Graham and G. Kolesnik. Van der Corput's Method of Exponential Sums. Cambridge University Press, New York, 1991.
\smallskip
\noindent
[Hu1990]  M. N.  Huxley, {\it Exponential Sums and Lattice Points.}
Proc. London Math. Soc. 60, 471-502, 1990 and
{\it Corrigenda: 'Exponential Sums and Lattice Points'.} Proc. London
Math. Soc. 66, 70, 1993.
\smallskip
\noindent
[Hu1997] M.N. Huxley,  Area lattice points and exponential sums,
London Mathematical Society Monographs, new series, 13, Oxford Science
Publications, The Clarendon Presss, Oxford University Press, New York, 1996
\smallskip
\noindent
[Hu2003] M.N. Huxley, {\it Exponential Sums and Lattice Points III,}
Proc. London Math. Soc. (2003).
\smallskip\noindent
[IM] H. Iwaniec and C.J.Mazzochi, {\it On the divisor and circle problems,}
J. Number Theory 29 (1988), no 1, 60-93
\smallskip
\noindent
[L] E. Landau. Vorlesung \"{u}ber Zahlentheorie, vols 1-3. Chelsea Publishing, 1947-50 [c. 1927], New York.
\smallskip\noindent

\smallskip\noindent
[P] G. Pick. {\it Geometrisches zur Zahlenlehre,}Sitzungber. Lotos Prag. 2, {\bf 19} (1870), 311-319.
\smallskip\noindent
[R] H. Rademacher. Topics in Analytic Number Theory, Springer-Verlag, Berlin, 1973.
\smallskip\noindent
[Sierp1] W. Sierpinski. {\it O pewnem zagadnieniu z rachunku funkcyj asymptotycznych } Prace mat. fyz. 17 (1906) pp. 77-118;
Summary in French: {\it Sur un probleme du calcul des fonctions
asymptotiques.} pp. 115-118.
\smallskip\noindent
[Sierp2] A. Schinzel. {\it Waclaw Sierpinski's Work in Number Theory} (Polish), Wiadomosci Matematyczne 26 (1) (1984), 24-31.
\smallskip\noindent
[T] E.C. Titchmarsh. The Theory of the Riemann Zeta-Function. Second Edition revised by
D.R. Heath-Brown. Clarendon Press, Oxford, 1986.
\smallskip\noindent
[VdC1] J.G. Van Der Corput. {\it Verscharfung der Abschatzungen beim Teilerproblem.}  Math. Ann. 87 (1922), 39-65
\smallskip\noindent
[VdC2] J.G. Van Der Corput. {\it Neue zahlentheoretische Abschatzungen.} Math. Ann. 89 (1923), pp. 215-254-65
\smallskip\noindent
[Vor] G. Voronoi. {\it Sur un probleme du calcul des fonctions asymptotiques.} J.r.a.M. 126 (1903) pp. 241-282
\smallskip\noindent
\bigskip
\noindent
Department of Mathematics, University of Pennsylvania

\noindent
209 South 33rd Street, Philadelphia, Pennsylvania, 19104

\end